\documentclass[12pt]{article}

\usepackage{amsmath}
\usepackage{amssymb}
\usepackage{amsthm}

\usepackage[dvips]{graphicx}
\usepackage{color}
\usepackage{epsfig}
\usepackage[mathscr]{eucal}
\newtheorem{theorem}{Theorem}

\newtheorem{prop}[theorem]{Proposition}

\newtheorem{lemma}[theorem]{Lemma}

\newcommand{\floor}[1]{\ensuremath{\left \lfloor {#1} \right \rfloor}}
\newcommand{\ceil}[1]{\ensuremath{\left \lceil {#1} \right \rceil}}

\newcommand{\cC}{{\cal C}}
\newcommand{\cD}{{\cal D}}
\newcommand{\cF}{{\cal F}}

\newcommand{\cI}{{\cal I}}
\newcommand{\cL}{{\cal L}}
\newcommand{\cP}{{\cal P}}

\newcommand{\cS}{{\cal S}}

\newcommand{\cX}{{\cal X}}
\newcommand{\cY}{{\cal Y}}

\newcommand{\bx}{{\bf x}}

\allowdisplaybreaks

\title{A Permutation Regularity Lemma}
\author{Joshua N. Cooper\footnote{Supported by NSF grant DMS-0303272.} \\ \small Courant Institute of Mathematics \\ \small New York University}
\date{\today}

\begin{document}

\maketitle

\begin{abstract}
We introduce a permutation analogue of the celebrated Szemer\'{e}di Regularity Lemma, and derive a number of consequences.  This tool allows us to provide a structural description of permutations which avoid a specified pattern, a result that permutations which scatter small intervals contain all possible patterns of a given size, a proof that every permutation avoiding a specified pattern has a nearly monotone linear-sized subset, and a ``thin deletion'' result.  We also show how one can count sub-patterns of a permutation with an integral, and relate our results to permutation quasirandomness in a manner analogous to the graph-theoretic setting.
\end{abstract}

\section{Introduction}

The Szemer\'{e}di Regularity Lemma, a tool developed in the early 1970's in service of the combinatorial milestone now known as the Szemer\'{e}di Theorem, has turned out to be one of the most useful tools in graph theory ever discovered.  In essence, it says that any graph can be approximated by a small collection of random-like graphs.  This powerful structural characterization allows one to answer questions about graphs by taking such a ``Szemer\'{e}di partition'' and then addressing the question by using known facts about random graphs.  A number of variants of the Regularity Lemma (or Uniformity Lemma, as it is sometimes called) have emerged since the publication of the original.  Versions of it giving structural decompositions of hypergraphs have been used in many contexts, and a few results have addressed the difficult case of sparse graph regularity.  The reader is encouraged to read the excellent surveys of Koml\'{o}s and Simonovits \cite{KS96} and Kohayakawa and R\"{o}dl \cite{KR03} to learn about how and where the Lemma is used, how it is proved, and what its limitations are.

An idea intimately related to regularity -- quasirandomness -- was introduced by Chung, Graham, and Wilson in \cite{CGW1}.  They show that a surprisingly large collection of random-like properties of graphs are in fact equivalent.  Then, in a series of remarkable papers, Chung and Graham applied similar analyses to hypergraphs, subsets of $\mathbb{Z}_n$, tournaments, and other combinatorial objects.  The following decade witnessed a flurry of generalizations and elaborations appearing in the literature, with much of the work exploring the connections between regularity and quasirandomness.  In particular, Simonovits and S\'{o}s \cite{SS91} showed how quasirandomness is equivalent to the property of having a Szemer\'{e}di partition into pieces whose regular pairs have density $1/2$.

The author defined quasirandom permutations in \cite{C04} and proved that several classes of simple arithmetic functions almost always give rise to quasirandom permutations (\cite{C03}).  The central paradigm is the same: a large collection of natural, random-like properties are mutually equivalent.  However, the connections with regularity break down in this realm, as it has not been possible so far to bridge the worlds of graph quasirandomness and permutation quasirandomness.

In the present paper, we remedy this situation by proving a regularity lemma for permutations and analogizing the basic results used alongside the graph Regularity Lemma.  The main result (Theorem \ref{uniformity}) says that the ground set of any permutation may be decomposed into a small exceptional set and a bounded number of intervals in the remaining points so that the action of the permutation is randomlike on each such interval.  Our hope is that this tool will help address the nascent realm of ``extremal permutation'' problems and lead to other work analogous to that of Extremal Graph Theory.\\

Examples of extremal permutation problems include:

\begin{enumerate}
\item For any permutation $\tau$, give a structural description of the permutation that {\it avoids} $\tau$, i.e., $\sigma |_I$ is not order-isomorphic to $\tau$ for any index set $I$.  The problem of showing that the number of such $\sigma$ is at most exponential in the number of symbols is commonly known as the ``Stanley-Wilf Conjecture'', and was recently solved by Marcus and Tardos \cite{MT04}.
\item For a given permutation $\tau$, which permutation $\sigma$ has the maximum number of ``copies'' of $\tau$, in the above sense?  We write the number of such copies as $\Lambda^\tau(\sigma)$.  This question has seen a number of advances in the past ten years, following Herb Wilf's address at the 1992 SIAM Conference on Discrete Math.  One particularly nice addition to the recent literature in this realm is \cite{HSV04}.
\item Given permutations $\tau$ and $\tau^\prime$, what is the expected value of $\Lambda^{\tau^\prime}(\sigma)$ in the space of permutations $\sigma$ chosen uniformly among those permutations on $n$ symbols which avoid $\tau$?  What is the maximum value of $\Lambda^{\tau^\prime}(\sigma)$ among all those permutations $\sigma$ which avoid $\tau$?
\item Call a sequence of permutations $\{\sigma_i\}_{i=1}^\infty$, $\sigma_i$ a permutation of $n_i$ symbols with $n_i \rightarrow \infty$, {\it asymptotically $k$-symmetric} if, for each $\tau$, a permutation on $k$ symbols, $\Lambda^\tau(\sigma_i) = \binom{n_i}{k}(1 + o(1))/k!$.  Does there exist, for all $k$, a sequence which is asymptotically $k$-symmetric but not asymptotically $(k+1)$-symmetric?  This question of R. L. Graham appears in \cite{C04} and is open except for $k=1,2,3$.
\end{enumerate}

The rest of the paper is as follows.  In the next section, we define regularity and uniformity for permutations and prove the existence of a regular/uniform partition.  Then, in Section \ref{patternavoidancesection}, we address Problem 1 above with structural results about permutations which avoid a given pattern.  These results are used in Section \ref{destroysection} to show that only a small number of pairs of points need be deleted to destroy all copies of a pattern in a permutation which has few of them to begin with.  Section \ref{quasirandomsection} provides a connection between permutation quasirandomness and regularity, and a proof of a new characterization of permutation quasirandomness.  The following section contains a discussion of the (asymptotic) pattern counts one can compute given a regular partition of any permutation, and the final section contains a full proof of the permutation regularity lemma. 

\section{Regularity} \label{regularitysection}

We provide two versions of a permutation regularity result, the latter of which appears to be the more interesting and applicable, and we distinguish the two settings through the use of the terms ``regular'' and ``uniform.''  The first result, concerning regularity, we state below but relegate the proof -- which is quite standard -- to Section \ref{proofsection}.

We consider permutations to be elements of $\mathfrak{S}_n$, the set of bijective maps from $\mathbb{Z}_n$ to itself.  For a permutation $\sigma \in \mathfrak{S}_n$ and subsets $S, T \subset \mathbb{Z}_n$, we write $p(S,T) = |\{(s,t)\in S \times T : \sigma(s) < t\}|$, and $d(S,T) = p(S,T)/|S||T|$.  Throughout the rest of this paper, we consider only partitions in which each $C_i$, $i \geq 1$, is an interval.  Though it is something of an abuse, we will often speak of a ``partition of $\sigma$'' instead of a partition of $\mathbb{Z}_n$.  For integers $s$ and $t$ and an $\epsilon > 0$, say that the pair $(C_s,C_t)$ is {\it $\epsilon$-regular} if, for all intervals $I \subset C_s$ and $J \subset C_t$ with $|I| \geq \epsilon |C_s|$ and $|J| \geq \epsilon |C_t|$, we have
$$
| d(I,J) - d(C_s,C_t) | \leq \epsilon.
$$
Then we call $\cP$ an {\it $\epsilon$-regular partition} into $k$ parts if $|C_s| = |C_t|$ for all $1 \leq s,t \leq k$, $|C_0| \leq \epsilon n$, and $(C_s,C_t)$ is an $\epsilon$-regular pair for all but $\epsilon k^2$ pairs $(s,t)$.  (If $P$ has only this first property, it is called {\it equitable}.)

Our first theorem is the following.

\begin{theorem}[Permutation Regularity] \label{regularity} Given $m \in \mathbb{N}^+$ and $\epsilon > 0$, there exist $M = M(\epsilon,m)$ and $N = N(\epsilon,m)$ so that any $\sigma : \mathbb{Z}_n \rightarrow \mathbb{Z}_n $ has an $\epsilon$-regular partition into $k$ (nonexceptional) intervals with $m \leq k \leq M$ if $n \geq N$.
\end{theorem}

Note that this statement is very similar to the one gotten by taking applying the ``standard'' Regularity Lemma for graphs to the bipartite graph whose color classes are two copies of $\mathbb{Z}_n$, and so that there is an edge from $s$ to $t$ if $\sigma(s) < t$.  The difference lies primarily in that the blocks of the partition must be intervals, and the two partitions of the color classes are actually the same. 

We now prove a reformulation of this result which will be easier to use for some applications.  Let $\cL(S,\alpha)$, for a set $S \subset \mathbb{Z}_n$ and $\alpha \in [0,1]$, denote the fraction of elements of $S$ whose image is less than $\alpha n$, i.e., $|\sigma(S) \cap [0,\alpha n)|/|S|$.  We say that two functions $f, g: [0,1] \rightarrow [0,1]$ are {\it $\epsilon$-near} if, for each $\alpha \in [0,1]$, $g(\alpha-\epsilon)-\epsilon \leq f(\alpha) \leq g(\alpha+\epsilon)+\epsilon$.  (We employ the convention that $g(\alpha)=g(0)$ for $\alpha < 0$ and $g(\alpha)=g(1)$ for $\alpha>1$.)  It is easy to see that this definition is symmetric in $f$ and $g$.

Now, we say that a partition $\{C_j\}_{j=0}^k$ of $\sigma$ is $(\epsilon,\cF)$-uniform, where $\cF = \{f_s\}_{s=1}^k$, if it is equitable, $|C_0| \leq \epsilon n$, and, for each $s \in \{1,\ldots,k\}$ and every interval $I \subset C_s$ with $|I| \geq \epsilon |C_s|$, $\cL(I,\cdot)$ is $\epsilon$-near $f_s$.

The following theorem, which we consider to be the main one of this paper, says essentially that permutations are, up to small deviations, concatenations of ``deterministic'' maps (ones which send all points into just a few small intervals) and ``random'' maps (ones which resemble the original map on each subinterval).  Note the absence of any ``exceptionality'' other than the exceptional set itself, in contrast to the Graph Regularity Lemma, where exceptional pairs are unavoidable.

\begin{theorem}[Permutation Uniformity] \label{uniformity} Given $m \in \mathbb{N}^+$ and $0 < \epsilon < 1$, there exists $M = M(\epsilon,m)$ and $N = N(\epsilon,m)$ so that, if $n \geq N$, $\sigma : \mathbb{Z}_n \rightarrow \mathbb{Z}_n $ has an $(\epsilon,\cF)$-uniform partition $\{C_j\}_{j=0}^k$, with $m \leq k \leq M$, where $\cF$ is a collection of $k$ nondecreasing $C^\infty$ functions $f_j : [0,1] \rightarrow [0,1]$.
\end{theorem}

\begin{proof} Without loss of generality, we may assume that $\epsilon < 1/2$.  Apply Theorem \ref{regularity}, and take an $\epsilon^2/4$-regular partition of $\sigma$ so that each $C_j$, $j \geq 1$, has cardinality $\leq \epsilon n/4$.  (We may always do so by choosing a partition with even higher regularity if necessary.)  Note that there can be at most $\epsilon k/2$ indices $s \in [k]$ so that there are more than $\epsilon k/2$ indices $t$ with $(s,t)$ being $\epsilon$-irregular.  Call all other $s$ ``good'', and add each ``bad'' $C_j$ to $C_0$ to create a new partition of $\sigma$.  Then the new exceptional set has size at most $\epsilon^2 n /4 + (\epsilon k / 2)(n/k) \leq \epsilon n$.

Fix a good $s$, and let $A$ be any subset of $\mathbb{Z}_n$.  Now, suppose $x, y \in \mathbb{Z}_n$ have the property that there is some $C_t \subset [x,y)$.  Then
$$
|C_t|^{-1} |\{(s,t) \in A \times C_t : \sigma(s) < t \}| \geq |\sigma(A) \cap [0,x)|,
$$
so $d(A,C_t) \geq \cL(A,x/n)$.  Similarly, $d(A,C_t) \leq \cL(A,y/n)$.  In order to guarantee that there is such a $C_t$ and that $(s,t)$ is regular, it suffices to ensure that the gap between $x$ and $y$ is at least
$$
(\epsilon k/2 + 1) |C_1| + \epsilon^2 n/4 \leq \epsilon n/2 + \epsilon n/4 + \epsilon^2 n/4 < \epsilon n,
$$
since it should be the length of $\epsilon k/2 + 1$ $C_j$'s plus all the points of $C_0$.  Therefore, if we set $y = (\alpha + \epsilon) n$, $x = \alpha n$, we have
$$
\cL(X,\alpha) \leq d(X,C_t) \leq \cL(X,\alpha+\epsilon).
$$
On the other hand, we may take $z = (\alpha - \epsilon)n$, and there will be a $C_{t^\prime} \subset [z,x)$, so that
\begin{equation} \label{squeeze}
\cL(X,\alpha-\epsilon) \leq d(X,C_{t^\prime}) \leq \cL(X,\alpha) \leq d(X,C_t) \leq \cL(X,\alpha+\epsilon).
\end{equation}
If we take $X=I$, an interval of $C_s$ of length at least $\epsilon |C_s|$, then we may apply (\ref{squeeze}) to get
$$
d(I,C_{t^\prime}) \leq \cL(I,\alpha) \leq d(I,C_t).
$$
Then, using the regularity of the partition, we see that
$$
d(C_s,C_{t^\prime}) - \epsilon \leq \cL(I,\alpha) \leq d(C_s,C_t) + \epsilon.
$$
Applying (\ref{squeeze}) once more, this time with $X = C_s$,
$$
\cL(C_s,\alpha-\epsilon) - \epsilon \leq \cL(I,\alpha) \leq \cL(C_s,\alpha+\epsilon) + \epsilon.
$$
Since this analysis works for any $\alpha \in [\epsilon,1-\epsilon)$, and the conclusion holds trivially otherwise, we may take $f_s(\alpha) = \cL(C_s,\alpha)$.\\

Note that $\cL(C_s,\alpha+1/n) - \cL(C_s,\alpha) \leq |C_s|^{-1} < 2k n^{-1}$.  It is easy to see, then, that by choosing $n$ large enough we may assume that all of the $f_s$ are $C^\infty$ and monotone.
\end{proof}

\section{Pattern Avoidance} \label{patternavoidancesection}

Define $\Lambda^\tau(\sigma)$ for $\tau \in \mathfrak{S}_m$ and $\sigma \in \mathfrak{S}_n$ to be the number of occurrences of the pattern $\tau$ in $\sigma$, i.e., the number of ``index sets'' $\{x_0 < \ldots < x_{m-1}\} \subset \mathbb{Z}_n$ such that $\sigma(x_i) < \sigma(x_j)$ iff $\tau(i) < \tau(j)$.

Suppose that $\sigma \in \mathfrak{S}_n$ has a uniform partition $P$, and $\tau \in \mathfrak{S}_m$.  If it is known that $\Lambda^\tau(\sigma)=o(n^m)$, what can be said about the $f_s$?  In fact, something quite strong: that it concentrates almost all the mass of $\sigma(C_s)$ in at most $m-1$ very small intervals.

\begin{theorem} \label{mblocks} Suppose $\sigma \in \mathfrak{S}_n$, $\tau \in \mathfrak{S}_m$, $0 < \epsilon \leq (2m)^{-1}$, and $n$ is sufficiently large.  Choose $\{C_j\}_{j=0}^k$, an $(\epsilon,\cF)$-uniform partition of $\sigma$.  If $\Lambda^\tau(\sigma) < (\epsilon n/2km)^m$, then, for each $1 \leq s \leq k$, there is a collection $\cI$ of at most $m-1$ disjoint intervals in $[0,1)$, each of length at most $6 \epsilon$, so that $|\sigma(C_s) \cap (n \cdot \bigcup \cI)| \geq |C_s| (1-7 m \epsilon)$.
\end{theorem}
\begin{proof} Write $\cF = \{f_s\}_{s=1}^k$.  First we prove a claim: if $J_0, \ldots, J_{m-1}$ are disjoint intervals of $[0,1)$ which are separated from each other by at least $4 \epsilon$, then, for some $t$, we have
$$
f_s(\sup J_t) - f_s(\inf J_t) \leq 5 \epsilon.
$$
To see this, suppose the contrary, i.e., that there are $m$ such intervals for which $f_s(\sup J_t) - f_s(\inf J_t) \geq 5 \epsilon$.  Then split $C_s$ into $m$ intervals $C_s^0, \ldots, C_s^{m-1}$ whose sizes differ by at most $1$, and denote their density functions by $f_s^q(\cdot) = \cL(C_s^q,\cdot)$.  Writing $x_t = \inf J_t$ and $y_t = \sup J_t$, we have
\begin{align*}
f_s^q(y_t + 2 \epsilon) - f_s^q(x_t - 2 \epsilon) & \geq (\cL(C_s,y_t + \epsilon)-\epsilon) - (\cL(C_s,x_t - \epsilon)+\epsilon) \\
& \geq (f_s(y_t) - \epsilon) - (f_s(x_t) + \epsilon) - 2\epsilon \\
& = f_s(y_t) - f_s(x_t) - 4 \epsilon \geq \epsilon,
\end{align*}
since $|C_s^q|/|C_s| \geq \epsilon$.  Define $x_t^\prime = \max\{0,x_t-2\epsilon\}$ and $y_t^\prime = \min\{1,y_t+2\epsilon\}$, and note that the intervals $\{J^\prime_t = [x_t^\prime, y_t^\prime)\}$ are disjoint, by the separation property of the $J_t$.  Then the fact that $f_s^q(y_t^\prime) - f_s^q(x_t^\prime) \geq \epsilon$ for each $q$ and $t$ implies that $|\sigma(C_s^q) \cap J^\prime_{\tau(q)}| \geq \epsilon |C_s^q| \geq \epsilon n/2km$.  If we take any $z_q \in C_s^q \cap \sigma^{-1}(J^\prime_{\tau(q)})$, then $z_0,\ldots, z_{m-1}$ is a $\tau$-pattern in $\sigma$, so we have at least $(\epsilon n/2km)^m$ such patterns, a contradiction.

Now, consider the following process: begin at $0$, and find the first $r$ so that $f_s(r) = 5 \epsilon$ (or $r=1$ if such a point does not exist).  This is possible because $f_s$ is monotone and continuous and $f_s(1)=1$.  Define $I_0 = [0,r)$.  Then, let $I^\prime_0 = [r,r+4 \epsilon)$.  Now, begin at $r+4\epsilon$, find the first $r^\prime$ so that $f_s(r^\prime) - f_s(r + 4\epsilon) = 5 \epsilon$ (or $r^\prime=1$, again, if this is not possible), and define $I_1 = [r+4\epsilon,r^\prime)$ and $I_1^\prime = [r^\prime, r^\prime+4 \epsilon)$.  Then define $r^{\prime \prime}$, $I_2$, and $I_2^\prime$ similarly, and so on.  This process must terminate in no more than $\ceil{1/(5 \epsilon)}$ steps, at which point the right-endpoint of the last interval defined is $1$.  In fact, it must terminate even sooner, by the claim above: if we have reached $I_m^\prime$, then $I_0, \ldots, I_{m-1}$ provide a contradiction.  Then the $I_l^\prime$ number at most $m-1$ and each has length at most $4 \epsilon$.  Now, define $I_l^{\prime\prime}$ to be the interval with left endpoint $x_l^{\prime\prime} = \min\{\inf I_l + \epsilon,1\}$ and right endpoint $y_l^{\prime\prime} = \max\{\sup I_l - \epsilon,0\}$.  Then
$$
\cL(C_s,y_l^{\prime\prime}) - \cL(C_s,x_l^{\prime\prime}) \leq (f_s(\sup I_l) + \epsilon) - (f_s(\inf I_l) - \epsilon) \leq 7 \epsilon.  
$$
Therefore, the intervals which comprise the complement of $\bigcup_l I_l^{\prime\prime}$, each of which contains some $I_l^\prime$, satisfy the conclusions of the theorem.
\end{proof}

Define a permutation $\sigma \in \mathfrak{S}_n$ to be {\it $m$-universal} if $\Lambda^\tau(\sigma) > 0$ for each $\tau \in \mathfrak{S}_m$.  Now, we say that a permutation $\sigma$ has the $(\delta, \epsilon, \gamma)$-property if, for every interval $I$ with $|I| \geq \delta n$ and every interval $J$ with $|J| \leq \epsilon n$, we have $|\sigma(I) \cap J| \leq \gamma |I|$.  That is, no sufficiently large interval is mapped too densely into any small interval.  Our next result says that, for the appropriate parameters, this property implies universality.  Note that, if we had instead stated that $|\sigma(I) \cap J| \geq \gamma |I|$ when $|J| \geq \epsilon n$, this would be immediate.  With the reverse inequalities, however, it is far from obvious.  On the other hand, if $\delta \gamma \geq \epsilon$, the statement would be vacuous.  Therefore, in particular, it has content whenever $\delta < \epsilon$.

\begin{prop} For each $m \geq 2$ and $\epsilon > 0$, there is a positive $\delta < \epsilon$ so that, for $n$ sufficiently large, if $\sigma \in \mathfrak{S}_n$ has the $(\delta,\epsilon,m^{-1})$-property, then $\sigma$ is $m$-universal.
\end{prop}
\begin{proof} Suppose the contrary, so that there is some $\tau \in \mathfrak{S}_m$ with $\Lambda^\tau(\sigma)=0$.  Take $\epsilon^\prime = \min\{\epsilon/6,m^{-2}/8\}$, and choose an $(\epsilon^\prime,\cF)$-uniform partition $\{C_j\}_{j=0}^k$.  Let $C_s$ be any block of the partition.  Then, by Theorem \ref{mblocks}, at least $|C_s| (1-7 m \epsilon^{\prime})/(m-1)$ points of $I = C_s$ are mapped by $\sigma$ into some interval $J$ of length at most $\epsilon n$.  However, if we take $\delta = 1/k$ and $\gamma = 1/(m-1)$, then the fact that $\sigma$ has the $(\delta,\epsilon,\gamma)$-property provides a contradiction, since $|I| \geq \delta n$, $J \leq \epsilon n$, and
$$
\frac{|\sigma(I) \cap J|}{|I|} \geq \frac{1-7m \epsilon}{m-1} > \frac{1}{m}.
$$
\end{proof}

Now, we show that any permutation which avoids a given $\tau$ has a linear sized subpattern which is ``nearly monotone''.  (Compare to the Erd\H{o}s-Szekeres Theorem, which says that any permutation on $n$ symbols has a $\sqrt{n}$-sized truly monotone subpattern.)  Define a permutation $\rho \in \mathfrak{S}_r$ to be {\it $\delta$-pseudomonotone} if either $\Lambda^{(01)}(\rho) \leq \delta \binom{r}{2}$ or $\Lambda^{(10)}(\rho) \leq \delta \binom{r}{2}$.  Then we have the following.

\begin{prop} \label{pseudomonotone} For every $\delta > 0$ and $\tau \in \mathfrak{S}_m$, there is a $c > 0$ so that, for any permutation $\sigma \in \mathfrak{S}_n$ which avoids $\tau$, with $n$ sufficiently large, there is a set $I \subset \mathbb{Z}_n$ with $|I| \geq cn$ so that $\sigma|_I$ is $\delta$-pseudomonotone.
\end{prop}
\begin{proof}  We may assume $m \geq 2$, and fix $\epsilon = \eta/14m$ with $\eta \leq 1$.  By Theorem \ref{mblocks}, $\sigma$ has an $(\epsilon,\cF)$-uniform partition so that, for each $1 \leq s \leq k$, there is an interval $I_s$ of length at most $6 \epsilon$ so that $|\sigma(C_s) \cap (n \cdot I_s)| \geq |C_s| (1-7 m \epsilon)/(m-1)$.  Order the $C_s$ left-to-right.  Suppose that $T$ of the $I_s$ intersect some fixed $I_t$.  At least
$$
T \cdot \frac{1-7m \epsilon}{m-1} \cdot |C_s| \geq \frac{T(1-7m \epsilon)(1-\epsilon)n}{k(m-1)}
$$
points are mapped by $\sigma$ into an interval of length at most $18 \epsilon n$.  Therefore,
$$
T \leq \frac{18 \epsilon k (m-1)}{(1-7m\epsilon)(1-\epsilon)} < 6 \eta k.
$$
Hence, we may iteratively pick $s_1, \ldots, s_{\ceil{\eta^{-1}/6}}$ so that the $I_{s_j}$ are mutually disjoint.  By the Erd\H{o}s-Szekeres Theorem, there is a subset $s^\prime_1 < \cdots < s^\prime_R$ of these $s_j$ of size at least $\ceil{\eta^{-1}/6}^{1/2}$ which is monotone with respect to the obvious ordering on the $I_{s_j}$.  Let $X$ be the union of the $\sigma^{-1}(I_{s^\prime_j}) \cap C_{s^\prime_j}$.  The only pairs of elements of $X$ which possibly display the opposite ordering to that of the intervals $I_{s^\prime_j}$ are contained within a single set of the form $\sigma^{-1}(I_{s^\prime_j}) \cap C_{s^\prime_j}$.  The fraction of pairs in $X$ of this type is at most
$$
\frac{R \binom{|C_s|}{2}}{\binom{R|C_s|/2(m-1)}{2}} \leq \frac{4(m-1)^2}{R} \leq 10 \eta^{1/2} (m-1)^2.
$$
If we let $\eta = \delta^2 (m-1)^{-4}/100$, then the set $X$ is $\delta$-pseudomonotone and has cardinality at least $n \cdot (20(m-1)/\delta k)$.
\end{proof}

\section{Destroying Patterns} \label{destroysection}

With the graph regularity lemma, one can prove that, if a graph $G$ contains at most $o(n^m)$ copies of some $m$-vertex graph, then we may remove $o(n^2)$ edges to destroy {\it all} copies.  Is there any hope of proving something analogous for permutations?

The first observation to make is that this is certainly not possible if one wishes to delete elements of the ground set.  Consider the permutation
$$
\sigma = (1,0,3,2,5,4,\ldots,n-1,n-2)
$$
for $n$ even.  It is clear that, even though this permutation has $\Lambda^{10}(\sigma) = O(n) = o(n^2)$, one must remove $\Omega(n)$ points to destroy all copies.  Furthermore, the generalization of this construction to other patterns is a simple matter.

Clearly, deleting points of the ground set is not the proper analogue of removing edges.  Let us instead attempt to ``delete'' pairs of points.  We wish to choose a subset $\cS \subset \binom{\mathbb{Z}_n}{2}$ so that {\it every copy of $\tau$ in $\sigma$ contains (in its index set) both points of some element of $\cS$}.  To state it another way: if we do not count index sets in which pairs from $\cS$ appear, there are no copies of the pattern $\tau$.  Any copy of $\tau$ containing such a pair we say is {\it destroyed} by the deletion of $\cS$.  The main result of this section says that, using $o(n^2)$ such deletions, we may destroy {\it all} copies of $\tau$ in a permutation $\sigma$ which has $\Lambda^\tau(\sigma) = o(n^m)$.

\begin{prop} Suppose that $\sigma \in \mathfrak{S}_n$, $\tau \in \mathfrak{S}_m$, and $\Lambda^{\tau}(\sigma) = o(n^m)$.  Then we may delete at most $o(n^2)$ index pairs to destroy all copies of $\tau$.
\end{prop}
\begin{proof} Take an $(\epsilon, \cF)$-uniform partition $\{C_j\}_{j=0}^k$, $\epsilon < (2m)^{-1}$, and choose $n$ large enough that $\Lambda^{\tau}(\sigma) < (\epsilon n/2km)^m$.  By Theorem \ref{mblocks}, for each $1 \leq s \leq k$, there is a collection $\cI_s$ of at most $m-1$ disjoint intervals $\{I_s^j\}$ in $[0,1)$, each of length at most $6 \epsilon$, so that $|\sigma(C_s) \cap (n \cdot \bigcup \cI_s)| \geq |C_s| (1-7 m \epsilon)$.  We create a new collection of families $\cI^\prime_s$ of intervals as follows.  Begin with the $\cI_s$.  If an interval $n \cdot \cI_s$ receives fewer than $\epsilon |C_s|$ points of $C_s$ under the action of $\sigma$, we remove it from the collection.  Then each $\cI^\prime_s$ has at most $m-1$ elements and $|\sigma(C_s) \cap (n \cdot \bigcup \cI^\prime_s)| \geq |C_s| (1-8 m \epsilon)$.

Now, delete all pairs which contain at least one point of $C_0$ or a point whose image does not fall into any of the $n \cdot \cI^\prime_s$.  There are most $(8m+1) \epsilon n^2$ of these.  Then, delete all pairs which contain two points from any one of the sets $C_s$.  This uses at most $k \binom{|C_s|}{2} \leq n^2/2k$ pairs.  Finally, delete all pairs whose elements are mapped to points within $12 \epsilon n$ of each other by $\sigma$.  There are at most $12 \epsilon n^2$ of these.  Hence, letting $\epsilon \rightarrow 0$ and $k \rightarrow \infty$, the result follows if we can show that the chosen deletion indeed destroys all copies of $\tau$.

Suppose not.  Then the index set on which $\tau$ appears, $i_0 < \ldots < i_{m-1}$ must have the following properties:
\begin{enumerate}
\item For each $r = 0, \ldots, m-1$, $\sigma(i_r) \in n \cdot I_s^j$ for some $s$ and $j$.
\item For each $q = 0, \ldots, m-2$, $\sigma(i_{\tau^{-1}(q)}) < \sigma(i_{\tau^{-1}(q+1)}) - 12 \epsilon n$.
\item If $\sigma(i_r) \in n \cdot I_s^j$ and $\sigma(i_{r^\prime}) \in n \cdot I_{s^\prime}^{j^\prime}$, then $s \neq s^\prime$.
\end{enumerate}
Since each of the $I_s^j$ have diameter at most $6 \epsilon$, the first two properties imply that the $n \cdot I_s^j$ must be disjoint.  Order these (dilated) intervals by increasing $s$, i.e., $s_0 < \cdots < s_{m-1}$, and call them $J_0,\ldots,J_{m-1}$.  Because they are disjoint and $\sigma(i_r) \in J_r$ for each $r$, the intervals themselves are ordered like a copy of $\tau$.  Therefore, since the indices $s$ are distinct, for {\it any} set of $m$ indices drawn one from each of $\sigma^{-1}(J_0) \cap C_{s_0}, \cdots, \sigma^{-1}(J_{m-1}) \cap C_{s_{m-1}}$, $\sigma$ restricted to this set is a copy of $\tau$.  This ensures that
$$
\Lambda^\tau(\sigma) \geq (\epsilon |C_s|)^m \geq (\epsilon n/2k)^m,
$$
a contradiction.
\end{proof}

\section{Quasirandomness} \label{quasirandomsection}
In \cite{C04}, the author proves that a number of random-like properties of permutations are equivalent to one another.  In order to state the main result of that paper, a few definitions are necessary.  Fix a permutation $\sigma \in \mathfrak{S}_n$.  For any $S, T \subset \mathbb{Z}_n$ we define the {\it discrepancy} of $S$ in $T$ as
$$
D_T(S) = \left | | S \cap T | - \frac{|S||T|}{n} \right |,
$$
and we define the discrepancy of a permutation $\sigma$ by
$$
D(\sigma) = \max_{I,J} D_J(\sigma(I)),
$$
where $I$ and $J$ vary over all intervals of $\mathbb{Z}_n$.  Also, define
$$
D^*(\sigma) = \max_{I,J} D_J(\sigma(I)),
$$
where $I$ and $J$ vary only over ``initial'' intervals, i.e., intervals of the form $[0,x)$.

We say that a sequence $\{\sigma_i\}_{i=1}^\infty$ of permutations of $\mathbb{Z}_{n_1}$, $\mathbb{Z}_{n_2}$, $\ldots$ is {\it quasirandom} if $D(\sigma_i) = o(n_i)$.  Often the indices are suppressed, and we simply say that $D(\sigma) = o(n)$.

By $e(x)$, we mean $e^{2 \pi i x}$.  We also use the convention that the name of a set and its characteristic function are the same.  The following is a portion of the main theorem in \cite{C04}.

\begin{theorem} For any sequence of permutations $\sigma \in \mathfrak{S}_n$, integer $m \geq 2$, and real $\alpha > 0$, the following are equivalent: \\[.1in]
{\bf [UB]} \parbox[t]{5.0in}{(Uniform Balance) $D(\sigma) = o(n)$.} \\[.13in]
{\bf [UB*]} \parbox[t]{5.0in}{(Uniform Star-Balance) $D^\ast(\sigma) = o(n)$.} \\[.13in]
{\bf [SP]} \parbox[t]{5.0in}{(Separability) For any intervals $I,J,K,K^\prime \subset
\mathbb{Z}_n$,
$$
\left | \sum_{x \in K \cap \sigma^{-1}(K^\prime)} I(x)
J(\sigma(x)) - \frac{1}{n} \sum_{x \in K,y \in K^\prime} I(x)J(y)
\right | = o(n)
$$} \\[.13in]
{\bf [mS]} \parbox[t]{5.0in}{(m-Subsequences) For any permutation $\tau \in \mathfrak{S}_m$ and intervals
$I,J \subset \mathbb{Z}_n$ with $|I| \geq n/2$ and $|J| \geq n/2$, we have $|I \cap \sigma^{-1}(J)| \geq n/4 + o(n)$ and
$$
\Lambda^\tau(\sigma |_{I \cap \sigma^{-1}(J)}) = \frac{1}{m!}
\binom{|\sigma(I) \cap J|}{m} + o(n^m).
$$} \\[.13in]
{\bf [2S]} \parbox[t]{5.0in}{(2-Subsequences) For any intervals $I,J \subset \mathbb{Z}_n$ with $|I| \geq n/2$ and $|J| \geq n/2$, we have $|I \cap \sigma^{-1}(J)| \geq n/4 + o(n)$ and
$$
\Lambda^{(01)}(\sigma |_{I \cap \sigma^{-1}(J)}) - \Lambda^{(10)}(\sigma |_{I \cap \sigma^{-1}(J)}) = o(n^2).
$$} \\[.13in]
{\bf [E($\alpha$)]} \parbox[t]{5.0in}{(Eigenvalue Bound $\alpha$) For all nonzero $k \in \mathbb{Z}_n$ and any interval $I$,
$$
\sum_{s \in \sigma(I)} e(-ks/n) = o(n|k|^{\alpha}).
$$} \\[.13in]
{\bf [T]} \parbox[t]{5.0in}{(Translation) For any intervals $I,J$,
$$
\sum_{k \in \mathbb{Z}_n} \left ( |\sigma(I) \cap (J+k)| - \frac{|I||J|}{n} \right )^2 = o(n^3).
$$}

Furthermore, for any implication between a pair of properties above, there exists a constant $K$ so that the error term $\epsilon_2 n^k$ of the consequent is bounded by the error term $\epsilon_1 n^l$ of the antecedent in the sense that $\epsilon_2 = O(\epsilon_1^K)$.
\end{theorem}

In \cite{SS91}, the authors connect graph quasirandomness and regularity by showing that, essentially, a sequence of graphs is quasirandom if and only if they possess density $1/2$ regular partitions with arbitrarily small $\epsilon$.  Here, we prove an analogous result for permutations.  Let $O_1(x)$ denote some real number whose absolute value is at most $x$.

\begin{prop} A sequence of permutations $\sigma_i \in \mathfrak{S}_{n_i}$, $i \geq 1$, $n_i \rightarrow \infty$, is quasirandom if and only if, for each $\epsilon > 0$, given any $(\epsilon,\{f_s\})$-uniform partition of $\sigma_i$ with $i$ sufficiently large, $f_s(x)$ is $2\epsilon$-near $id(x)=x$ for each $s$.
\end{prop}
\begin{proof} Suppose that $\sigma$ is quasirandom, and let $P$ be an $(\epsilon,\cF)$-uniform partition.  For any interval $C_s$, {\bf [UB]} implies that
$$
\cL(C_s,\alpha) = \alpha + o(1).
$$
Therefore, if we choose $n_i$ large enough, we may ensure that $|\cL(C_s,\alpha) - \alpha| \leq \epsilon$ for all $\alpha \in [0,1]$, which immediately implies that $f_s$ is $2\epsilon$-near $id$ for each $s$.

On the other hand, suppose $\sigma_i$ has an $(\epsilon,\{f_s\})$-uniform partition $P = \{C_j\}_{j=0}^k$ for all sufficiently large $i$, where $f_s$ is $2\epsilon$-near $id$ for each $s$.  It is easy to see that this implies (by sub-additivity) that $P$ is $(3\epsilon,\{id\})$-uniform.  We may assume that the $C_j$ are ordered left-to-right.  Choose $[0,x], [0,y) \subset \mathbb{Z}_{n_i}$.  If, for some $s$, $C_j \subset [0,x)$ for all $1 \leq j \leq s$, choose the largest such $s$, and let $X = [0,x) \cap C_{s+1}$ (or $\emptyset$ if $s=k$).  Otherwise, let $X = [0,x)$ and $s=0$.  Then
$$
[0,x) = \left (\bigcup_{1 \leq j \leq s} C_j \right ) \cup X \cup E
$$
for some $E \subset C_0$.  Therefore, $s|C_1| + |X| = x + O_1(\epsilon n_i)$.  We may write
$$
\bigg |\sigma \big ( [0,x) \big) \cap [0,y) \bigg| = \left |\bigcup_{1 \leq j \leq s} \sigma(C_j)  \cap [0,y) \right | + |\sigma(X) \cap [0,y)| + |\sigma(E) \cap [0,y)|.
$$
Clearly, $\big |\sigma(E) \cap [0,y) \big| = O_1(\epsilon n_i)$.  Define $\beta = y/n_i$.  Applying uniformity, then, 
\begin{align*}
\bigg |\sigma \big ([0,x) \big) \cap [0,y) \bigg| & = \sum_{1 \leq j \leq s} |C_1| \cL(C_j,\beta) + |X| \cL(X,\beta) + O_1(\epsilon n_i)\\
& = \sum_{1 \leq j \leq s} |C_1| (\beta + O_1(3 \epsilon)) + |X| (\beta + O_1(3 \epsilon)) + O_1(\epsilon n_i) \\
& = (s|C_1| + |X|) \beta + O_1(4 \epsilon n_i),
\end{align*}
so we may conclude that
$$
\left| \big |\sigma \big ([0,x) \big) \cap [0,y) \big | - \frac{xy}{n_i} \right| \leq 5 \epsilon n_i.
$$
If we take $\epsilon \rightarrow 0$, then $\sigma_i$ is quasirandom, by {\bf [UB*]}.

\end{proof}

\section{Counting Subpatterns} \label{patterncountssection}

We wish to count how many occurrences of the pattern $\tau \in \mathfrak{S}_m$ appear in a permutation $\sigma$ of $\mathbb{Z}_n$ with a given $(\epsilon,\cF)$-uniform partition $\{C_j\}_{j=0}^k$.  Unfortunately, we have no control over the structure of the exceptional set $C_0$, so it is not possible to get an ``exact'' count this way.  Nonetheless, if we write $X = \bigcup_{j} C_j$, it is easy to see that, for $n \geq \epsilon^{-1} \geq 2$,
\begin{equation} \label{eq1}
\left | \binom{n}{m} - \binom{|X|}{m} \right | \leq 2 \epsilon n^m/(m-1)!
\end{equation}
so that a count of the $\tau$-patterns on $X$ is going to be close to the same count on all of $\mathbb{Z}_n$.  Note that we may also ignore all but the set of occurrences of $\tau$ all of whose symbols occur in different $C_j$'s, since the number of these is $\binom{k}{m} |C_s|^m$, which is off from $\binom{n}{m}$ by at most
\begin{equation} \label{eq2}
\left | \binom{n}{m} - \binom{k}{m} |C_s|^m \right | \leq 2 \epsilon n^m/(m-1)! 
\end{equation}
for $n$ sufficiently large.

There is an additional obstruction to counting patterns that is more subtle than these two issues.  Suppose the mass of $\sigma(C_s)$, for some $s$, were very tightly concentrated in some interval.  If we use $f_s$ as an estimate of its density function, then, since the condition of $\epsilon$-nearness can ``dislocate'' the entire mass of $\sigma(C_s)$ by up to $\epsilon$, the counts could be off by a significant amount.  On the other hand, only index sets whose images have two points close to one another can be affected in this way.  Since there are few of these, with some work, we are able to ignore them in the total count.

One way ensure that the counts are accurate is simply to posit that $f_s$ does not concentrate its mass too tightly.  Therefore, define $f_s$ to be $(B,\epsilon)$-Lipschitz if, for each $x$, $|f_s(x + \epsilon) - f_s(x)| \leq B \epsilon$.  For example, if we have a (quasi-)random permutation, we may take $f_s(x) = x$ for each $s$, a function which is $(1,\epsilon)$-Lipschitz for each $\epsilon$.

The following lemma makes this idea rigorous.  Define
$$
\mathfrak{D}^r(\beta)=\{\bx=(x_1,\ldots,x_r) \in [0,\beta)^r : \forall i, 1 \leq i < r, x_i < x_{i+1}\}.
$$

\begin{lemma} \label{integralestimate} Let $f_j, g_j$, $j = 1,\ldots,r$, be cumulative distribution functions on $[0,a]$, $a > 0$, and suppose that for each $j$, $f_j$ is $\epsilon$-near $g_j$.  If $g_j$ is $(B,\epsilon)$-Lipschitz for each $j$, then
$$
\left | \int_{\mathfrak{D}^r(\beta)} \alpha(x_1) \, dF_r - \int_{\mathfrak{D}^r(\beta)} \alpha(x_1) \, dG_r \right | \leq r(a+1)(B+1)\epsilon.
$$
for any $0 \leq \beta \leq a$, where $dF_r = df_1(x_1) \cdots df_r(x_r)$, $dG_r = dg_1(x_1) \cdots dg_r(x_r)$, and $\alpha:[0,a] \rightarrow [0,1]$ is any nondecreasing function.
\end{lemma}
\begin{proof} We prove the result by induction, by repeated application of integration by parts.  First, we check that it is true for $r=1$.  Note that, since $f_j$ is $\epsilon$-near $g_j$,
$$
f_j(\beta) = g_1(\beta) + O_1((B+1) \epsilon)
$$
by the Lipschitz property, and we have
\begin{align*}
\int_{\mathfrak{D}^1(\beta)} \alpha(x_1) df_1(x_1) &= \int_0^\beta \alpha(x_1) df_1(x_1) \\
& = \alpha(\beta) f_1(\beta) - \int_0^\beta f_1(x_1) d\alpha(x_1) \\
& = \alpha(\beta) g_1(\beta) + O_1((B+1) \epsilon) \\
& \qquad - \int_0^\beta [g_1(x_1) + O_1((B+1)\epsilon)] d\alpha(x_1) \\
& = \int_0^\beta \alpha(x_1) dg_1(x_1) + O_1((a+1)(B+1)\epsilon).
\end{align*}
Now, suppose the result holds for $r-1$, with $r > 1$.  Then
\begin{align*}
\int_{\mathfrak{D}^r(\beta)} dF_r &= \int_0^\beta \int_{\mathfrak{D}^{r-1}(x_{r})} \alpha(x_1) \, dF_{r-1} \, df_r(x_r) \\
&= \int_0^\beta \left ( \int_{\mathfrak{D}^{r-1}(x_{r})} \alpha(x_1) dG_{r-1} + O_1((r-1)(a+1)(B+1)\epsilon) \right ) \, df_r(x_r) \\
&= \int_0^\beta \int_{\mathfrak{D}^{r-1}(x_{r})} \alpha(x_1) \, dG_{r-1} \, df_r(x_r) + O_1((r-1)(a+1)(B+1)\epsilon)
\end{align*}
The function $\alpha_1(x_r) = \int_{\mathfrak{D}^{r-1}(x_{r})} \alpha(x_1) \, dG_{r-1}$ is nondecreasing, nonnegative, and bounded by $1$, so we may apply the $r=1$ case to get
$$
\int_{\mathfrak{D}^r(\beta)} dF_r = \int_{\mathfrak{D}^r(\beta)} dG_r + O_1(r(a+1)(B+1)\epsilon).
$$
\end{proof}

We wish to be able to count subpatterns in permutations which do not necessarily have the Lipschitz property, however.  In order to be able to use this result, we have the following Lemma which says that convolving the c.d.f. of a permutation with a uniform distribution on a short interval preserves nearness and gives us a Lipschitz property.  Therefore, fix $\delta >0$ and, given a c.d.f. $f$ on $[0,1]$, define $\tilde{f}(t) = \delta^{-1} \int_{-\delta}^0 f(t+s) ds$, a c.d.f. on $[0,1+\delta]$.

\begin{lemma} \label{convolve} If $f$ and $g$, c.d.f.'s on $[0,1]$, are $\epsilon$-near, then $\tilde{f}$ and $\tilde{g}$ are $\epsilon$-near.  Furthermore, $\tilde{f}$ is $(2 \epsilon \delta^{-1}, \epsilon)$-Lipschitz.
\end{lemma}
\begin{proof} To see the first claim, we write
$$
\tilde{g}(t+\epsilon)-\tilde{f}(t) = \delta^{-1} \int_{t-\delta}^t g(s+\epsilon) - f(s) \, ds \leq \epsilon.
$$
For the second claim,
\begin{align*}
\tilde{f}(t+\epsilon)-\tilde{f}(t) &= \delta^{-1} \left ( \int_{t+\epsilon-\delta}^{t+\epsilon} f(s) \,ds - \int_{t-\delta}^{t} f(s) \,ds \right) \\
&=\delta^{-1} \left ( \int_{t}^{t+\epsilon} f(s) \,ds - \int_{t-\delta}^{t+\epsilon-\delta} f(s) \,ds \right) \\
& \leq 2 \epsilon \delta^{-1}.
\end{align*}
\end{proof}

Now, fix a permutation $\tau \in \mathfrak{S}_m$.  Write $\mathfrak{D}$ for $\mathfrak{D}^m(1)$, $\mathfrak{D}^\prime$ for $\mathfrak{D}^m(1+\delta)$, and define a differential form $d\omega_f$ on $[0,1)^m$ as follows:
$$
d\omega_f = |C_1|^m \sum_{1 \leq s_0 < \ldots < s_{m-1} \leq k} \bigwedge_{j=0}^{m-1} df_{s_{\tau^{-1}(j)}}(x_j)
$$
We define $d\omega_g$, etc., analogously.  Then we have the following.

\begin{theorem} \label{count} Suppose $d\omega_f$ is defined as above, $\epsilon \leq 1/2$, and $n$ is sufficiently large.  Then
$$
\left |\Lambda^\tau(\sigma) - \int_\mathfrak{D} d\omega_f \right| < (20 \epsilon^{1/2} m^2 + 4/k) n^m/(m-1)!
$$
\end{theorem}
\begin{proof} Define $g_j(x) = \cL(C_j,x)$.  Suppose that $m$ elements $x_i \in \mathbb{Z}_n$ are chosen uniformly at random, $x_i \in C_{j_i}$, for some sequence $j_0 < \ldots < j_{m-1}$.  Writing $j^\prime_r = j_{\tau^{-1}(r)}$, the probability that their images under $\sigma$ form a $\tau$ is precisely
$$
\int_\mathfrak{D} dg_{j^\prime_0}(x_1) \cdots dg_{j^\prime_{m-1}}(x_m)
$$
because $dg_j$ represents the distribution of the images $\sigma(x_i)$ for $x_i$ chosen at random from $C_j$, and we wish to compute the probability that
$$
\sigma(x_{j^\prime_0}) < \cdots < \sigma(x_{j^\prime_{m-1}}).
$$
Multiplying by the number $|C_1|^m$ of $m$-tuples, adding over all subsets $\{j_i\}$, and accounting for (\ref{eq1}) and (\ref{eq2}) yields
\begin{equation} \label{eq3}
|\Lambda^\tau(\sigma) - \int_\mathfrak{D} d\omega_g| \leq 4\epsilon n^m/(m-1)!
\end{equation}
By Lemma \ref{integralestimate} and Lemma \ref{convolve} with $\delta = \epsilon^{1/2}$, the quantity
$$
\left | \int_{\mathfrak{D}^\prime} d\tilde{g}_{j^\prime_0}(x_1) \cdots d\tilde{g}_{j^\prime_{m-1}}(x_m) \! - \! \int_{\mathfrak{D}^\prime} d\tilde{f}_{j^\prime_0}(x_1) \cdots d\tilde{f}_{j^\prime_{m-1}}(x_m) \right |
$$
is bounded by $m(2+\epsilon^{1/2})(2\epsilon^{1/2}+1)\epsilon \leq 9m \epsilon$.  Summing up again, we find
\begin{equation} \label{eq4}
\left|\int_{\mathfrak{D}^\prime} d\omega_{\tilde{f}} - \int_{\mathfrak{D}^\prime} d\omega_{\tilde{g}} \right| \leq 9\epsilon n^m/(m-1)!
\end{equation}
Now, let $dG = k^{-1} \sum_{j=1}^k dg_j$ and $d\tilde{G} = k^{-1} \sum_{j=1}^k d\tilde{g}_j$.  If we choose an $m$-tuple of points from each of these distributions, the distributions of their orderings with respect to increasing $j$ coincide so long as each point is at least $\delta$ away from all the others, since we may view $d\tilde{G}$ as a random draw from the distribution $dG$ followed by a random ``jump'' forward uniformly distributed in $[0,\delta]$.  The probability that such an $m$-tuple has two points at most $\delta$ apart is bounded by the probability that some pair of its points are that close, i.e.,
\begin{align*}
& \leq \binom{m}{2} \max_x \int_{x-\delta}^{x+\delta} dG \\
& = \binom{m}{2} \max_x \frac{1}{k} \sum_{j=1}^k \int_{x-\delta}^{x+\delta} dg_j \\
& = \binom{m}{2} \max_x \frac{1}{k} \sum_{j=1}^k \left ( \cL(C_j,x+\delta)-\cL(C_j,x-\delta) \right ) \\
& = \binom{m}{2} \max_x \frac{|\bigcup_j \sigma(C_j) \cap [(x-\delta)n,(x+\delta)n)|}{k|C_1|} \\
& \leq \frac{m^2}{2} \cdot \frac{2 \delta n}{k |C_1|} \leq m^2 \delta (1-\epsilon)^{-1} \leq 2 m^2 \delta.
\end{align*}
Since $d\omega_{g}$ (or $d\omega_{\tilde{g}}$) is the same as the distribution of an unordered $m$-tuple drawn from $dG$ (resp., $d\tilde{G}$) {\it minus} the event that two points are drawn from the same $C_j$,
\begin{equation} \label{eq5}
\left|\int_\mathfrak{D} d\omega_{g} - \int_{\mathfrak{D}^\prime} d\omega_{\tilde{g}} \right| \leq (2 \epsilon^{1/2} m^2 + 2/k) n^m/m!,
\end{equation}
where the second summand follows from the fact that
$$
1 - \binom{k}{m} \cdot \frac{m!}{k^m} \leq 1 - (1 - \frac{m-1}{k})^m < 1/k.
$$
Similarly, if we define $dF = k^{-1} \sum_{j=1}^k df_j$, $d\tilde{F} = k^{-1} \sum_{j=1}^s d\tilde{f}_j$, $F(x) = \int_0^x dF$, and $G(x) = \int_0^x dG$, the difference of the two integrals in question is bounded by $\binom{m}{2} \binom{n}{m}$ times the probability that two points chosen from $dF$ are within $\delta$, or
\begin{align*}
& \leq \binom{m}{2} \max_x \int_{x-\delta}^{x+\delta} dF \\
&= \binom{m}{2} \max_x (F(x+\delta)-F(x-\delta)) \\
& \leq \binom{m}{2} (\max_x (G(x+\delta+\epsilon)-G(x-\delta-\epsilon))+2\epsilon) \\
& \leq m^2 (\delta+2\epsilon) (1-\epsilon)^{-1} \leq 6 \epsilon^{1/2} m^2 
\end{align*}
And so,
\begin{equation} \label{eq6}
\left|\int_\mathfrak{D} d\omega_{f} - \int_\mathfrak{D^\prime} d\omega_{\tilde{f}} \right| \leq (6 \epsilon^{1/2} m^2 + 2/k) n^m/m!
\end{equation}
Putting together (\ref{eq3}), (\ref{eq4}), (\ref{eq5}), and (\ref{eq6}), we have
$$
\left |\Lambda^\tau(\sigma) - \int_\mathfrak{D} d\omega_f \right| < (20 \epsilon^{1/2} m^2 + 4/k) n^m/m!
$$
\end{proof}

\section{The Proof of Theorem \ref{regularity}} \label{proofsection}

For the proof of Theorem \ref{regularity} below, we are heavily indebted to \cite{D00}, which we find to have the most comprehensible -- if not the shortest -- proof of the Regularity Lemma in the literature.

For disjoint sets $X,Y \subset \mathbb{Z}_n$, define the ``index''
$$
q(X,Y)=|X||Y|d^2(X,Y)/n^2.
$$
Then extend this definition to a pair of partitions $\cX$ of $X$ and $\cY$ of $Y$ by $q(\cX,\cY) = \sum_{X^\prime \in \cX,Y^\prime \in \cY} q(X^\prime,Y^\prime)$.  For a partition $\cP = \{C_j\}_{j=1}^k$ of $\mathbb{Z}_n$, we write $q(\cP) = \sum_{i, j} q(C_i,C_j)$.  If one set in the partition, $C_0$, has been designated as an {\it exceptional set}, then we treat $C_0$ as a collection of singletons in this sum.  That is, we write $\tilde{\cP}$ for the partition which refines $\cP$ by splitting $C_0$ into singletons, then $q(\cP) = q(\tilde{\cP},\tilde{\cP})$.  First of all, note that
$$
q(\cP) = \sum_{i,j} \frac{|C_i||C_j|}{n^2} d^2(C_i,C_j) \leq n^{-2} \sum_{i,j} |C_i| |C_j| = 1.
$$

Now, we have the following simple lemma, which says that refinement can only increase the index of a partition.

\begin{lemma} \label{nodecrease}
\mbox{}
\begin{enumerate}
\item Let $C, D \subset \mathbb{Z}_n$ (not necessarily disjoint).  If $\cC$ is a partition of $C$ and $\cD$ is a partition of $D$, then $q(\cC,\cD) \geq q(C,D)$.
\item If $\cP,\cP^\prime$ are partitions of $\mathbb{Z}_n$ and $\cP^\prime$ refines $\cP$, then $q(\cP^\prime) \geq q(\cP)$.
\end{enumerate}
\end{lemma}
\begin{proof}
\mbox{}
\begin{enumerate}
\item Let $\cC = \{C_1,\ldots,C_k\}$ and $\cD = \{D_1,\ldots,D_l\}$.  Then
\begin{align*}
q(\cC,\cD) &= \frac{1}{n^2} \sum_{i,j} \frac{p(C_i,D_j)^2}{|C_i||D_j|} \\
& \geq \frac{1}{n^2} \frac{( \sum_{i,j} p(C_i,D_j) )^2}{\sum_{i,j} |C_i||D_j|} \\
& = \frac{1}{n^2} \frac{p(C,D)^2}{(\sum_i |C_i|)(\sum_j |D_j|)} \\
& = q(C,D),
\end{align*}
where the inequality follows from Cauchy-Schwarz.
\item Let $\cP = \{C_1,\ldots,C_k\}$, and for $1 \leq i \leq k$, let $\cC_i$ be the partition of $C_i$ induced by $\cP^\prime$.  Then
$$
q(\cP) = \sum_{i, j} q(C_i,C_j) \leq \sum_{i, j} q(\cC_i,\cC_j) = q(\cP^\prime),
$$
where the inequality follows from part (1).
\end{enumerate}
\end{proof}

The next lemma says that we may exploit irregular pairs to increase the index somewhat.

\begin{lemma} Let $\epsilon > 0$, and let $C,D \subset \mathbb{Z}_n$ be intervals.  If $(C,D)$ is not $\epsilon$-regular with respect to $\sigma$, then there are partitions $\cC$ and $\cD$ of $C$ and $D$, respectively, so that
$$
q(\cC,\cD) \geq q(C,D) + \epsilon^4 \frac{|C||D|}{n^2}.
$$
\end{lemma}
\begin{proof} Suppose $(C,D)$ is not $\epsilon$-regular, and choose intervals $C_1 \subset C$ and $D_1 \subset D$ with $|C_1| > \epsilon |C|$ and $|D_1| > \epsilon |D|$ so that $|\eta| > \epsilon$, where $\eta = d(C_1,D_1)-d(C,D)$.  Let $\cC = \{C_1,C_2,C_3\}$ and $\cD = \{D_1,D_2,D_3\}$, where $C_2$ is the ``left half'' of $C \setminus C_1$; $C_3$ is the ``right half''; and $D_2$ and $D_3$ are defined similarly.  (That is, $C_1$ splits the interval $C$ into three pieces: $C_1$ itself, one interval of points less than those of $C_1$ and one interval of points greater than those of $C_1$.  Either, but not both, of these may be empty.)

For ease of notation, write $c_i = |C_i|$, $d_i = |D_i|$, $e_{ij} = p(C_i,D_j)$, $c = |C|$, $d = |D|$, and $e = p(C,D)$.  Then, applying Cauchy-Schwarz again, we see
\begin{align*}
q(\cC,\cD) &= \frac{1}{n^2} \sum_{i,j} \frac{e_{ij}^2}{c_id_j} = \frac{1}{n^2}  \left ( \frac{e_{11}^2}{c_1d_1} + \sum_{(i,j) \neq (1,1)} \frac{e_{ij}^2}{c_id_j} \right) \\
& \geq \frac{1}{n^2} \left ( \frac{e_{11}^2}{c_1d_1} + \frac{(e-e_{11})^2}{cd - c_1d_1} \right).
\end{align*}
Since $e_{11} = c_1d_1e/cd + \eta c_1d_1$, we have
\begin{align*}
n^2 q(\cC,\cD) &\geq \frac{1}{c_1d_1} \left (\frac{c_1d_1e}{cd}+\eta c_1d_1 \right)^2 + \frac{1}{cd - c_1d_1}\left(\frac{cd-c_1d_1}{cd}e-\eta c_1 d_1\right)^2 \\
&= \frac{e^2}{cd} + \frac{\eta^2 c_1 d_1 cd}{cd - c_1d_1} \\
&\geq \frac{e^2}{cd} + \epsilon^4cd
\end{align*}
since $c_1 \geq \epsilon c$, $d_1 \geq \epsilon d$, and $\eta^2 > \epsilon^2$. 
\end{proof}

The following lemma is the crux of the proof of Theorem \ref{regularity}.

\begin{lemma} \label{refinementstep} Let $0 < \epsilon \leq 1/4$ and $k \in \mathbb{N}$, let $\sigma$ be a permutation of $\mathbb{Z}_n$, and let $\cP$ be an equitable partition of $\mathbb{Z}_n$ into $\{C_j\}_{j=0}^k$ with $|C_0| \leq \epsilon n$ and $|C_j| \geq 81^k$ for $j > 0$.  If $\cP$ is not $\epsilon$-regular, then there is an equitable partition $\cP^\prime = \{C_j^\prime\}_{j=0}^l$ of $\mathbb{Z}_n$ with exceptional set $C_0^\prime$, where $k \leq l \leq k 81^k$, such that $|C_0^\prime| \leq |C_0| + n/9^k$ and
$$
q(\cP^\prime) \geq q(\cP) + \epsilon^5/2.
$$
\end{lemma}

\begin{proof} Let $c = |C_1|$.  For all $1 \leq i,j \leq k$, define a partition $\cC_{ij}^1$ of $C_i$ and a partition $\cC_{ji}^2$ of $C_j$ as follows.  If the pair $(C_i,C_j)$ is $\epsilon$-regular, then let $\cC^1_{ij} = \{C_i\}$ and $\cC^2_{ji} = \{C_j\}$.  If not, then by the previous lemma, there are tripartitions $\cC^1_{ij}$ and $\cC^2_{ji}$ of $C_i$ and $C_j$, respectively, so that
$$
q(\cC^1_{ij},\cC^2_{ji}) \geq q(C_i,C_j) +  \frac{\epsilon^4c^2}{n^2}.
$$
For each $i = 1, \ldots, k$, let $\cC_i$ be the partition of $C_i$ that is the common refinement of every partition $\cC^r_{ij}$.  Note that $|\cC_i| \leq 9^k$.  Now, consider the partition
$$
\cC = \{C_0\} \cup \bigcup_{i=1}^k \cC_i,
$$
with $C_0$ as exceptional set.  Then $\cC$ refines $\cP$ and $k \leq |\cC | \leq k 9^k$.

Let $\cC_0 = \{\{x\}:x \in C_0\}$.  If $\cP$ is not $\epsilon$-regular, then for more than $\epsilon k^2$ of the pairs $(C_i,C_j)$, the partitions $\cC^1_{ij}$ and $\cC^2_{ji}$ are nontrivial.  Hence, by Lemma \ref{nodecrease},
\begin{align*}
q(\cC) &= \sum_{i,j \geq 1} q(\cC_i,\cC_j) + \sum_{i \geq 0} q(\cC_0,\cC_i) \\
& \geq \sum_{i,j \geq 1} q(\cC^1_{ij},\cC^2_{ji}) + \sum_{i \geq 1} q(\cC_0,\{C_i\}) + q(\cC_0) \\
& \geq \sum_{i,j \geq 1} q(C_i,C_j) + \epsilon k^2 \frac{\epsilon^4 c^2}{n^2} + \sum_{i \geq 1} q(\cC_0,\{C_i\}) + q(\cC_0) \\
&= q(\cP) + \epsilon^5 \left ( \frac{kc}{n} \right )^2 \\
& \geq q(\cP) + \epsilon^5/2,
\end{align*}
since $kc \geq 3n/4$.

Now, $\cC$ satisfies the conclusions of the theorem, except that it may not be equitable.  To fix the situation, cut each non-exceptional block of $\cC$ into a maximal collection of disjoint intervals of size $d = \floor{c/81^k}$.  Call the resulting set of intervals $\{C_j^\prime\}_{j=1}^l$, and let $C_0^\prime = \mathbb{Z}_n \setminus \bigcup C_j^\prime$.  This new partition $\cP^\prime$ refines $\cC$, so
$$
q(\cP^\prime) \geq q(\cC) \geq q(\cP)+\epsilon^5/2.
$$
Since each set $C_j^\prime$, $j > 0$, is contained in one of the sets $C_i$, but not more than $81^k$ sets can lie inside the same $C_j$, we also have $k \leq l \leq k 81^k$.  On the other hand, the sets $C_1^\prime,\ldots,C_l^\prime$ use all but at most $d$ points from each nonexceptional block of $\cC$.  Therefore,
\begin{align*}
|C_0^\prime| &\leq |C_0| + d |\cC| \\
& \leq |C_0| + \frac{c}{81^k} k 9^k \\
& = |C_0| + ck 9^{-k} \\
& \leq |C_0| + n 9^{-k}.
\end{align*}
\end{proof}

Now, since $q(\cP) \leq 1$, this lemma cannot be applied ad infinitum.  Indeed, we may now complete the proof of Theorem \ref{regularity}.

\begin{proof}[Proof of Theorem \ref{regularity}]  Let $\epsilon > 0$ and $m \geq 1$.  Without loss of generality, $\epsilon \leq 1/4$.  Let $s = \ceil{2/\epsilon^5}$.  If we apply Lemma \ref{refinementstep} $s$ times, we end up with a partition of $\sigma$ which is $\epsilon$-regular.

However, we need to choose the ``starting'' partition so that the exceptional set ends up with cardinality $\leq \epsilon n$ and the $C_j$ are sufficiently large at each stage.  With each iteration of the lemma, the size of the exceptional set can grow by at most $n/9^k$.  Therefore, we wish to choose $k$ large enough so that $s$ increments of $n/9^k$ add up to at most $\epsilon n/2$, and $n$ large enough so that $|C_0| < k$ implies $|C_0| \leq \epsilon n/2$.  (We can guarantee $|C_0|<k$ if we begin with an equitable $k$-partition.)  So let $k \geq m$ be large enough so that $9^k \geq 2s/\epsilon$.  Then $s/9^k \leq \epsilon/2$, and hence
$$
k + \frac{s}{9^k}n \leq \epsilon n
$$
whenever $n \geq 2k/\epsilon$.

Now, define $f(x) = x81^x$.  We may take $M = \max\{f^s(k),2k/\epsilon\}$.  To deal with the second condition -- that the blocks be sufficiently large at each stage -- note that, after $s$ steps, the nonexceptional blocks sizes are at least $n/(2M)^s$.  Therefore, choosing $N = \max\{2M/\epsilon,81^M (2M)^s\}$ suffices, and the proof is complete.
\end{proof}

\section{Concluding Remarks}

The discussion of Section \ref{patternavoidancesection} is largely ``local'', i.e., the analysis is concerned with the internal structure of individual blocks of the uniform partition.  Section \ref{patterncountssection} consists of a ``global'' analysis -- it does not take into account the internal structure of the blocks, only their relationships with one another.  On the other hand, the proofs of the main results of Sections \ref{destroysection} and \ref{quasirandomsection}, as well as that of Proposition \ref{pseudomonotone}, are both.  It is here, in the interplay between local and global, that we believe the most interesting behavior resides.  We suspect that such dual analysis may lead to a better understanding of extremal permutations in the senses of Problem 2 and 3 of the Introduction, perhaps using the results of Section \ref{patterncountssection}.  Theorem \ref{count}, in theory, gives a translation of these problems from combinatorial to analytic.  We are hopeful that Theorem \ref{count} can find application in algorithmic settings, e.g., in the vein of \cite{DLR95}, or in other contexts where a ``counting lemma'' has been useful, such as the hypergraph-theoretic proof of the Szemer\'edi Theorem.

\section{Acknowledgements}

Thank you to Jim Propp, Vera S\'{o}s, and Joel Spencer for stimulating questions and invaluable discussions.  Thanks also to the referee for helpful comments and suggestions.


\begin{thebibliography}{100}
\bibitem{AF1} N.\@ Alon and E.\@ Friedgut, On the number of permutations avoiding a given pattern, {\it J. Comb. Theory Ser. A} {\bf 89} (2000), 133--140.
\bibitem{CGW1} F.\@ R.\@ K.\@ Chung, R.\@ L.\@ Graham, and R.\@ M.\@ Wilson, Quasi-random graphs, {\it Combinatorica} {\bf 9} (1989), 345--362.
\bibitem{C04} J.\@ N.\@ Cooper, Quasirandom permutations, {\it J. Comb. Theory Ser. A} {\bf 106} 1 (2004), 123--143.
\bibitem{C03} J.\@ N.\@ Cooper, Quasirandom arithmetic permutations,  J. Number Theory  {\bf 114}  (2005),  no. 1, 153--169. 
\bibitem{D00} R.\@ Diestel, Graph Theory. Second Edition. {\it Graduate Texts in Mathematics} {\bf 173}, Springer-Verlag, New York, 2000.
\bibitem{DLR95} R.\@ A.\@ Duke, H.\@ Lefmann, V.\@ R\"odl, A fast approximation algorithm for computing the frequencies of subgraphs in a given graph, SIAM J. Comput.  24  (1995),  no. 3, 598--620. 
\bibitem{HSV04} M.\@ Hildebrand, B.\@ E.\@ Sagan, and V.\@ R.\@ Vatter, Bounding quantities related to the packing density of $1(l+1)l\ldots 2$, preprint, 2004.
\bibitem{KR03} Y.\@ Kohayakawa, V.\@ R\"{o}dl, Szemer\'{e}di's regularity lemma and quasi-randomness,  {\it Recent advances in algorithms and combinatorics}, 289-–351, CMS Books Math./Ouvrages Math. SMC, {\bf 11}, Springer, New York, 2003.
\bibitem{KS96} J.\@ Koml\'{o}s and M.\@ Simonovits,  Szemer\'{e}di's regularity lemma and its applications in graph theory,  {\it Combinatorics -- Paul Erd\H{o}s is Eighty}, vol. 2, D. Mikl\'{o}s, V. T. S\'{o}s, and T.\@ Sz\H{o}nyi, eds., Bolyai Mathematical Studies, pages 295--352.  J\'{a}nos Bolyai Mathematical Society, Budapest, Budapest, 1996.
\bibitem{MT04} A.\@ Marcus and G.\@ Tardos, Excluded permutation matrices and the Stanley-Wilf conjecture, J. Combin. Theory Ser. A {\bf 107} (2004), no. 1, 153--160.
\bibitem{SS91} M.\@ Simonovits and V.\@ T.\@ S\'{o}s, Szemer\'{e}di's partition and quasirandomness, {\it Random Structures Algorithms} {\bf 2} (1991), no. 1, 1--10.
\bibitem{Sz76} E.\@ Szemer\'{e}di,  Regular partitions of graphs, {\it Probl\`{e}mes Combinatoires et Th\'{e}orie des Graphes}, Colloques Internationaux CNRS n. 260, Orsay, 1976.
\bibitem{W99} H.\@ S.\@ Wilf, The patterns of permutations, {\it Discrete Math.} {\bf 257} (2002), no. 2-3, 575–-583.
\end{thebibliography}
\end{document}